\documentclass[11pt,reqno]{article}
\usepackage{latexsym, amsmath, amssymb, amsthm, a4, epsfig}
\usepackage{graphicx}
\usepackage{mathrsfs}

\newtheorem{theorem}{Theorem}[section]
\newtheorem{cor}[theorem]{Corollary}
\newtheorem{lemma}[theorem]{Lemma}
\newtheorem{prop}[theorem]{Proposition}

\newcommand{\nm}{\noalign{\smallskip}}
\newcommand{\ds}{\displaystyle}

\newcommand{\p}{\partial}

\newcommand{\eqnref}[1]{(\ref {#1})}

\newcommand{\Cbb}{\mathbb{C}}

\newcommand{\Kbb}{\mathbb{K}}
\newcommand{\Rbb}{\mathbb{R}}
\newcommand{\Sbb}{\mathbb{S}}

\newcommand{\la}{\langle}
\newcommand{\ra}{\rangle}

\newcommand{\Hcal}{\mathcal{H}}
\newcommand{\Lcal}{\mathcal{L}}

\newcommand{\Hscr}{\mathscr{H}}
\newcommand{\Vscr}{\mathscr{V}}

\def\Ba{{\bf a}}
\def\Bb{{\bf b}}

\def\Bf{{\bf f}}

\def\Bn{{\bf n}}

\def\Bu{{\bf u}}
\def\Bv{{\bf v}}

\def\Bx{{\bf x}}
\def\By{{\bf y}}
\def\Bz{{\bf z}}
\def\BA{{\bf A}}
\def\BB{{\bf B}}

\def\BE{{\bf E}}
\def\BF{{\bf F}}
\def\BG{{\bf G}}

\def\BK{{\bf K}}

\def\BS{{\bf S}}

\newcommand{\Ga}{\alpha}

\newcommand{\Gd}{\delta}

\newcommand{\Gf}{\phi}
\newcommand{\Gvf}{\varphi}

\newcommand{\Gk}{\kappa}

\newcommand{\Gl}{\lambda}

\newcommand{\Gm}{\mu}
\newcommand{\Gv}{\nu}

\newcommand{\Gr}{\rho}

\newcommand{\Gs}{\sigma}

\newcommand{\Go}{\omega}

\newcommand{\Gy}{\psi}

\newcommand{\GF}{\Phi}
\newcommand{\GG}{\Gamma}

\newcommand{\GO}{\Omega}

\newcommand{\GY}{\Psi}

\newcommand{\BGG}{{\bf \GG}}

\newcommand{\BGvf}{\mbox{\boldmath $\Gvf$}}
\newcommand{\Bpsi}{\mbox{\boldmath $\Gy$}}
\newcommand{\BGF}{{\bf \GF}}

\newcommand{\beq}{\begin{equation}}
\newcommand{\eeq}{\end{equation}}

\def\ol{\overline}
\newcommand{\hatna}{\widehat{\nabla}}

\numberwithin{equation}{section}
\numberwithin{figure}{section}

\begin{document}

\title{Cloaking by anomalous localized resonance for linear elasticity on a coated structure\thanks{This work is supported by A3 Foresight Program among China (NSF), Japan (JSPS), and Korea (NRF 2014K2A2A6000567). Work of HK is supported by NRF 2016R1A2B4011304}}

\author{Kazunori Ando\thanks{Department of Electrical and Electronic Engineering and Computer Science, Ehime University, Ehime 790-8577, Japan. (ando@cs.ehime-u.ac.jp)} \and Hyeonbae Kang\thanks{Department of Mathematics, Inha University, Incheon
    22212, S. Korea (hbkang@inha.ac.kr)} \and Kyoungsun Kim\thanks{Department of Mathematical Sciences, Seoul National University, Seoul 08826, S. Korea (kgsunsis@snu.ac.kr)} \and Sanghyeon Yu\thanks{Seminar for Applied Mathematics, ETH Z\"urich, R\"amistrasse 101, CH-8092 Z\"urich, Switzerland (sanghyeon.yu@sam.math.ethz.ch)}}

\date{}

\maketitle

\begin{abstract}
We investigate anomalous localized resonance on the circular coated structure and cloaking related to it in the context of elasto-static systems. The structure consists of the circular core with constant Lam\'e parameters and the circular shell of negative Lam\'e parameters proportional to those of the core. We show that the eigenvalues of the Neumann-Poincar\'e operator corresponding to the structure converges to certain non-zero numbers determined by Lam\'e parameters and derive precise asymptotics of the convergence. We then show with estimates that cloaking by anomalous localized resonance takes place if and only if the dipole type source lies inside critical radii determined by the radii of the core and the shell.
\end{abstract}

\section{Introduction} \label{sec:intro}

If a dielectric material is coated by a meta-material with a negative dielectric constant, then a source outside the structure may cause an anomalous localized resonance by which the source may be cloaked. This phenomenon was first discovered in \cite{MN_PRSA_06, NMM_94} when the core and shell are concentric disks and the source is single or multiple polarizable dipoles. It was shown in \cite{MN_PRSA_06} that there is a critical radius determined by radii of the core and the shell such that if the dipole source is located inside the radius, then cloaking by anomalous localized resonance (CALR) takes place, and if the dipole source is located outside the radius, then CALR does not take place.

In a recent work \cite{ACKLM}, quantitatively precise analysis of CALR is presented on circular coated structure with arbitrary sources. Among the findings of the paper is that CALR is a spectral phenomenon occurring at the limit point of eigenvalues of the Neumann-Poincar\'e (NP) operator defined on concentric disks. In fact, by Plemelj's symmetrization principle, the NP operator can be realized as a symmetric operator \cite{KPS}. Since it is compact on concentric circles, the NP operator has real eigenvalues converging to $0$ (exponentially fast), and there CALR takes place. On a single disk $0$ is the only eigenvalue of the NP operator (other than $1/2$ of multiplicity $1$). So CALR does not take place. But by coating the disk with another disk, the boundary becomes two circles and the NP eigenvalues are perturbed and converges to $0$, and hence CALR does take place on the coated structure. In fact, CALR may take place on a single inclusion (not a coated inclusion). For example, the NP eigenvalues on ellipses are not $0$ and converge to $0$, CALR takes place on ellipses \cite{AK-JMAA-16}. There is yet another approach to CALR using the variational method, for which we refer to \cite{KLSW-CMP-14}.

Since CALR is a spectral phenomenon at the limit point of eigenvalues of the NP operator, it is interesting to look into CALR in the context other than electro-statics where the notion of the NP operator is defined. In this regard the system of elasto-statics is particularly interesting since the elastic NP operator is not compact even on smooth boundaries. However, it is proved in the recent paper \cite{AJKKY} that elastic NP eigenvalues on smooth boundaries of two-dimensional domains, whose Lam\'e parameters are $(\Gl, \Gm)$, accumulate at either $k_0$ or $-k_0$ where
\beq\label{kzero}
k_0= \frac{\Gm}{2 \left( \Gl + 2 \Gm \right)}.
\eeq

Using this spectral property of the elastic NP operator CALR for the linear isotropic elasticity is considered in \cite{AJKKY} when the background is the usual isotropic elastic material with the pair of Lam\'e constants $(\Gl, \Gm)$ and the inclusion has
\beq\label{meta}
  (\tilde{\Gl}, \tilde{\Gm}) = (c + i \Gd) (\Gl, \Gm)
\eeq
as its Lam\'e constants, where $c$ is a negative constant and $\Gd$ is a small positive constant representing dissipation. The asymptotics of NP eigenvalues are derived and it is proved that CALR occurs on ellipses if the constant $c$ in \eqnref{meta} satisfies
\beq\label{ccond}
z(c):= \frac{c+1}{2(1-c)} = \pm k_0.
\eeq
We emphasize that the condition \eqnref{ccond} can be satisfied only when $c$ is negative since $0<k_0 <1/2$. The elastic NP eigenvalues on two dimensional disks are
\beq\label{NPeigendisk}
  \frac{1}{2}, \quad - \frac{\Gl}{2 (2 \Gm + \Gl)}, \quad \pm k_0,
\eeq
where the first two have multiplicity one, while the last two are of infinite multiplicities, and CALR does not occur since $\pm k_0$ are the elastic NP eigenvalues (not limit points). One may expect that like the electro-static case if we coat the disk by an concentric circle, then the NP eigenvalues in \eqnref{NPeigendisk} are perturbed from $\pm k_0$ (while converging to them) and CALR may occur.

The purpose of this paper is to confirm it. To do so, we suppose that the disk $D_i$ is included in another concentric disk $D_e$ where the Lam\'e parameters are given by
\beq\label{total_lame_const}
  (\Gl^\Gd, \Gm^\Gd) =
  \begin{cases}
    (\Gl, \Gm), & \text{ in } \Rbb^2 \setminus \ol{D_e}, \\
    (c + i \delta) (\Gl, \Gm), & \text{ in } D_e \setminus \ol{D_i}, \\
    (\Gl, \Gm), & \text{ in } D_i.
  \end{cases}
\eeq
Here $c$ is a negative constant and $\Gd$ is a parameter which will be sent to $0$. We show that the elastic NP eigenvalues on the boundary of $D_e \setminus \ol{D_i}$ deviate from $\pm k_0$, but accumulate at $\pm k_0$. We actually show that the NP operator (with respect to the inner product on $H^{-1/2}$-space defined by the single layer potential) can be realized as a series of block matrices, and then using perturbation theory of eigenvalues we derive asymptotic formula of the convergence. We then show that CALR occurs if $c$ satisfies \eqnref{ccond}. It is worth emphasizing that the critical radii (for cloaking) when $z(c)=k_0$ and when $z(c)=-k_0$ are different, which is due to the different behavior of NP eigenvalues near $k_0$ and $-k_0$.

A few comments are in order before completing Introduction. The major emphasis of this paper is on the spectral property of the elastic NP operator on two circles. Once we obtain the asymptotic behavior of eigenvalues, the rest of arguments are parallel to those in \cite{ACKLM} even though they are more involved since a system of equations are being dealt with. Existence of elastic meta-materials satisfying \eqnref{meta} is out of the scope of our expertise. In this regard, we mention that the method of this paper works only when \eqnref{meta} is satisfied. If, for example, only one of the parameters is negative, the method may not work.

While this work is in progress, we were informed by Hongyu Liu that their work \cite{LL_SIMA-16} was completed. There it is proved CALR occurs on the annulus structure in elasticity when the source function $\Bf$ (see \eqnref{negative_indx_elasticity_problem}) is supported in a circle. We emphasize that the method of this paper is completely different from that of \cite{LL_SIMA-16}: their method uses the variational one similar to that of \cite{KLSW-CMP-14} while our method uses spectral properties of the NP operator. The method of this paper reveals more clearly the spectral nature of anomalous localized resonance.

This paper is organized as follows. Section \ref{sec:prelim} is to formulate the problem using layer potentials and define the elastic NP operator. Section \ref{sec:NP} is to derive asymptotics of the NP eigenvalues near $\pm k_0$. Occurrence of CALR is proved in section \ref{sec:CALR}. Appendix is to provide proofs of two technical identities.

\section{Layer potential formulation of the problem} \label{sec:prelim}

\subsection{Layer potentials for the Lam\'e system}\label{sec:layer}

Let us first recall definitions of the layer potential and the NP operator related to the Lam\'e system (see, for example, \cite{AK07}). With the pair of Lam\'{e} constants $(\Gl, \Gm)$ satisfying $\Gm > 0$ and $\Gl + \Gm > 0$, the isotropic elasticity tensor $\Cbb = \left( C_{ijkl} \right)_{i, j, k, l = 1}^2$ is defined by
\beq
  C_{ijkl} := \Gl \Gd_{ij} \Gd_{kl} + \Gm \left( \Gd_{ik} \Gd_{jl} + \Gd_{il} \Gd_{jk} \right). \nonumber
\eeq
Then the corresponding Lam\'{e} system of elasticity equations is defined to be $\Lcal_{\Gl, \Gm}:= \nabla \cdot \Cbb \hatna$
where $\hatna$ is the symmetric gradient, namely,
$$
\hatna \Bu := \frac{1}{2} \left( \nabla \Bu + \nabla \Bu^T \right) \quad (T \mbox{ for transpose}).
$$

Let $\GO$ be a connected bounded domain with the Lipschitz boundary in $\Rbb^2$. The single layer potential of the density function $\BGvf$ on $\p \GO$ associated with the Lam\'{e} operator $\Lcal_{\Gl, \Gm}$ is defined by
\begin{align*}
  \BS_{\p \GO} [\BGvf] (\Bx) & := \int_{\p \GO} \BGG(\Bx - \By) \BGvf(\By) d\Gs(\By), \quad \Bx \in \Rbb^2,
\end{align*}
where $d\Gs$ is the arc length of $\p \GO$ and $\BGG(\Bx) = \left( \GG_{ij}(\Bx) \right)_{i, j = 1}^2$ is the Kelvin matrix of the fundamental solution to the Lam\'{e} system in $\Rbb^2$, namely,
\beq
  \GG_{ij}(\Bx) = \frac{\Ga_1}{2 \pi} \Gd_{ij} \ln{\left\vert \Bx \right\vert} - \frac{\Ga_2}{2 \pi} \frac{x_i x_j}{\left\vert \Bx \right\vert^2},
\eeq
with
\beq
  \Ga_1 = \frac{1}{2} \left( \frac{1}{\Gm} + \frac{1}{2 \Gm + \Gl} \right) \ \text{ and } \ \Ga_2 = \frac{1}{2} \left( \frac{1}{\Gm} - \frac{1}{2 \Gm + \Gl} \right). \nonumber
\eeq

The (elasto-static) NP operator on $\p\GO$ is defined by
\beq
  \BK_{\p \GO}^* [\BGvf] (\Bx) := \mbox{p.v.} \int_{\p \GO} \p_{\Gv_\Bx} \BGG(\Bx - \By) \BGvf(\By) d\Gs(\By) \quad \text{a.e. } \Bx \in \p \GO. \nonumber
\eeq
Here $\mbox{p.v.}$ stands for the Cauchy principal value, the conormal derivative on $\p \GO$ associated with $\nabla \cdot \Cbb \hatna$ is defined to be
\beq
  \p_\Gv \Bu := ( \Cbb \hatna \Bu ) \Bn = \Gl \left( \nabla \cdot \Bu \right) \Bn + 2 \Gm ( \hatna \Bu ) \Bn \quad \text{ on } \p \GO,
\eeq
where $\Bn$ is the outward unit normal to $\p \GO$, and $\p_{\Gv_\Bx} \BGG(\Bx - \By)$ is defined by
\beq
  \p_{\Gv_\Bx} \BGG(\Bx - \By) \Bb = \p_{\Gv_\Bx} \left( \BGG(\Bx - \By) \Bb \right)
\eeq
for any constant vector $\Bb \in \Rbb^2$.

Let $H^{1 / 2}(\p \GO)^2$ be the usual $L^2$-Sobolev space of order $1 / 2$ and $H^{- 1 / 2}(\p \GO)^2$ be its dual space with respect to $L^2$-pairing $\la \cdot, \cdot \ra$.
Let $\GY$ be the subspace of $H^{- 1 / 2}(\p \GO)^2$ spanned by
\beq\label{Psi}
\begin{bmatrix} 1 \\ 0 \end{bmatrix},  \quad \begin{bmatrix} 0 \\ 1 \end{bmatrix}, \quad \begin{bmatrix} y \\ - x\end{bmatrix}.
\eeq
We denote by $H^{-1/2}_\GY(\p \GO)$ the annihilator of $\GY$, namely, the collection of all $\BGvf \in H^{-1/2}(\p \GO)^2$ such that $\la \BGvf, \Bf \ra = 0$ for all $\Bf \in \GY$.

The NP operator $\BK_{\p \GO}^*$ is bounded on $H^{-1 / 2}(\p \GO)^2$, and maps $H^{-1 / 2}_\GY(\p \GO)$ into itself. If $\Gl \notin [-1/2, 1/2)$, then $\Gl I - \BK_{\p \GO}^*$ is invertible on $H^{-1 / 2}(\p \GO)^2$ (and on $H^{-1 / 2}_\GY(\p \GO)$). The NP operator is related to the single layer potential through the following jump formula (see \cite{AK07, DKV88}):
\begin{align}
  \left. \p_\Gv \BS_{\p \GO} [\BGvf] \right|_\pm & = \left( \pm \frac{1}{2} I + \BK_{\p \GO}^* \right) [\BGvf] \quad \text{a.e. on } \p \GO, \label{single_jump_formula}
\end{align}
where the subscripts $\pm$ denote the limits (to $\p\GO$) from outside and inside of $\p \GO$. Even if $\BK_{\p \GO}^*$ is not self-adjoint with respect to the usual $L^2$-inner product, it can be realized as a self-adjoint operator on $H^{-1 / 2}_\GY(\p \GO)$. In fact, it is proved in \cite{AJKKY} (following the discovery of \cite{KPS}) that the product $(\cdot, \cdot)_*$, defined by
\beq
  \left( \BGvf, \Bpsi \right)_{*} := - \la \BGvf, {\BS}_{\p\GO} [\Bpsi] \ra,
\eeq
is actually an inner product on $H^{-1 / 2}_\GY(\p \GO)$ which yields a norm equivalent to the usual $H^{-1/2}$-norm. Then $\BK_{\p \GO}^*$ is  self-adjoint with respect to this new inner product thanks to the Plemelj's symmetrization principle:
$$
\BS_{\p\GO} \BK_{\p\GO}^* = \BK_{\p\GO} \BS_{\p\GO},
$$
where $\BK_{\p \GO}$ is the $L^2$-adjoint operator of $\BK_{\p \GO}^*$. It is well-known (see, for example, \cite{DKV88}) that $\BK_{\p \GO}^*$ is not a compact operator on $H^{-1 / 2}(\p \GO)$ even if $\p\GO$ is smooth. However, it is proved in \cite{AJKKY} that the spectrum of $\BK_{\p \GO}^*$ consists of pure point spectrum converging to $\pm k_0$, which is defined in \eqnref{kzero}.

\subsection{The CALR problem on coated disks}

To formulate the problem on coated disks, let $D_i := \{ |\Bx| < r_i \}$ and $D_e := \{ |\Bx| < r_e\}$, $0 < r_i < r_e$, and write $\GG_i=\p D_i$ and $\GG_e=\p D_e$. The distribution of Lam\'e parameters are given by \eqnref{total_lame_const}, and the elasticity tensor $\Cbb^\Gd = ( C_{ijkl}^\Gd )_{i,j,k,l=1}^2$ is given by
\beq
  C_{ijkl}^\Gd := \Gl^\Gd \Gd_{ij} \Gd_{kl} + \Gm^\Gd \left( \Gd_{ik} \Gd_{jl} + \Gd_{il} \Gd_{jk} \right).
\eeq
We then consider the following elasticity equation:
\beq\label{negative_indx_elasticity_problem}
    \nabla \cdot \Cbb^{\Gd} \hatna \Bu^{\Gd} = \Bf \quad \mbox{ in } \Rbb^2
\eeq
with the decaying condition $\Bu^\Gd(\Bx) \to {\bf 0}$ as $|\Bx| \to \infty$, where $\Bu^\Gd = ( u_1^\Gd, u_2^\Gd)^T$ and the source $\Bf$ is compactly supported in $\Rbb^2 \setminus \ol{D_e}$ which satisfies the following conservation law:
\beq
  \int_{\Rbb^2} \Bf = {\bf 0}. \label{conservation}
\eeq

Let $\Bu^\Gd$ be the solution to \eqnref{negative_indx_elasticity_problem} and define
\beq
  E(\Bu^\Gd):= \int_{D_e \setminus D_i} \Gl | \nabla \cdot \Bu^\Gd |^2 + 2 \Gm |\hatna \Bu^\Gd |^2.
\eeq
Here and afterwards, $\BA : \BB$ for two matrices $\BA = \left( a_{ij} \right)$ and $\BB = \left( b_{ij} \right)$ denotes $\sum_{i, j} a_{ij} b_{ij}$, and $|\BA|^2=\BA:\BA$. Then CALR is characterized by the following two conditions:
\begin{itemize}
\item the blow up of the dissipated energy on the annulus:
\beq\label{blow_up_diss_energy}
    E^\Gd := \Gd E(\Bu^\Gd) \to \infty \quad \text{ as } \Gd \to 0,
\eeq
\item the boundedness of the solution $\Bu^\Gd(\Bx)$ far away from the structure; more precisely, there is a radius $a > r_e$ such that $\vert \Bu^\Gd(\Bx) \vert < C$ for some $C > 0$ on $\vert \Bx \vert > a$ as $\Gd \to 0$.
\end{itemize}
It is worth mentioning that $E^\Gd$ is the imaginary part of the total energy, namely,
$$
E^\Gd = \Im \int_{\Rbb^2} \Cbb^\Gd \hatna\Bu^\Gd : \hatna\Bu^\Gd .
$$

\subsection{Layer potential formulation of the CALR problem}

We now formulate the CALR problem \eqnref{negative_indx_elasticity_problem} using layer potentials and define the NP operator for the problem. The formulation here is parallel to that in \cite{ACKLM} for the electro-static case. So we omit proofs.

Let us fix notation first: let $\Hscr^* = H^{-1 / 2}(\GG_i)^2 \times H^{-1 / 2}(\GG_e)^2$ and  $\Hscr^*_\GY:= H^{-1/2}_\GY(\GG_i)^2 \times H^{-1/2}_\GY (\GG_e)^2$. Let $\BF$ be the Newtonian potential of $\Bf$:
\beq
  \BF(\Bx) := \int_{\Rbb^2} \BGG(\Bx - \By) \Bf(\By) d\By, \quad \Bx \in \Rbb^2. \label{Fdef}
\eeq
Note that $\BF$ satisfies $\Lcal_{\Gl, \Gm} \BF = \Bf$ in $\Rbb^2$ and $\BF(\Bx) \to {\bf 0}$ as $\vert \Bx \vert \to \infty$ since $\Bf$ satisfies \eqref{conservation}. Let $\BS_{\GG_i}$ and $\BS_{\GG_e}$ be the single layer potentials on $\GG_i$ and $\GG_e$, respectively, with respect to the Lam\'e parameters $(\Gl, \Gm)$.
Then we seek the solution $\Bu^\Gd$ of \eqnref{negative_indx_elasticity_problem} in the following form:
\beq\label{rep_sol}
  \Bu^\delta(\Bx) = \BF(\Bx) + \BS_{\GG_i} [\BGvf_i^\Gd] (\Bx) + \BS_{\GG_e} [\BGvf_e^\Gd] (\Bx)
\eeq
for some $(\BGvf_i^\Gd ,\BGvf_e^\Gd) \in \Hscr^*_\GY$. We emphasize that the solution can be represented as \eqnref{rep_sol} because the Lam\'e parameters of the shell is of the form $(c + i \Gd) (\Gl, \Gm)$. If they are not of the form, then the single layer potential with different Lam\'e parameter should be used, and arguments of this paper may not be applied.

The transmission conditions to be satisfied by $\Bu^\Gd$ on $\GG_i$ and $\GG_e$ are
\beq
  \begin{cases}
    \left( c + i \Gd \right) \left. \p_\Gv \Bu^\Gd \right|_+ =\left. \p_\Gv \Bu^\Gd \right|_- & \text{ on } \GG_i, \\
    \left. \p_\Gv \Bu^\Gd \right|_+ = \left( c + i \Gd \right) \left. \p_\Gv \Bu^\Gd \right|_- & \text{ on } \GG_e.
  \end{cases} \nonumber
\eeq
So $(\BGvf_i^\Gd, \BGvf_e^\Gd)$ is the solution to the following system of integral equations:
\beq
  \begin{cases}
    & \left( c + i \Gd \right) \left. \p_{\Gv_i} \BS_{\GG_i} [\BGvf_i^\Gd] \right|_+ - \left. \p_{\Gv_i} \BS_{\GG_i} [\BGvf_i^\Gd] \right|_- + \left(c - 1 + i \Gd \right) \p_{\Gv_i} \BS_{\GG_e} [\BGvf_e^\Gd] \\
    &\qquad\qquad\qquad = \left(1 - c - i \Gd \right) \p_{\Gv_i} \BF \quad \text{ on } \GG_i, \\
    & \left(1 - c - i \Gd \right) \p_{\Gv_e} \BS_{\GG_i} [\BGvf_i^\Gd] + \left. \p_{\Gv_e} \BS_{\GG_e} [\BGvf_e^\Gd] \right|_+ - \left( c + i \Gd \right) \left. \p_{\Gv_e} \BS_{\GG_e} [\BGvf_e^\Gd] \right|_- \\
    &\qquad\qquad\qquad = \left(c - 1 + i \Gd \right) \p_{\Gv_e} \BF  \quad \text{ on } \GG_e,
  \end{cases} \nonumber
\eeq
where $\p_{\Gv_i}$ and $\p_{\Gv_e}$ are the conormal derivatives on $\GG_i$ and $\GG_e$, respectively.
Using the jump formula \eqref{single_jump_formula}, the above integral equations are written as
\beq
  \begin{bmatrix}
    - z_\Gd I + \BK_{\GG_i}^* & \p_{\Gv_i} \BS_{\GG_e} \\
    \p_{\Gv_e} \BS_{\GG_i} & z_\Gd I + \BK_{\GG_e}^*
  \end{bmatrix}
  \begin{bmatrix}
    \BGvf_i^\Gd \\
    \BGvf_e^\Gd
  \end{bmatrix}
  = -
  \begin{bmatrix}
    \p_{\Gv_i} \BF \\
    \p_{\Gv_e} \BF
  \end{bmatrix}, \label{raw_matrix_equation}
\eeq
where
\beq\label{zd}
  z_\Gd = \frac{1 + c + i \Gd}{2 \left( 1 - c \right) - 2 i \Gd}.
\eeq
We emphasize that $\p_{\Gv_i} \BF \in H^{-1/2}_\Psi(\GG_i)$ and $\p_{\Gv_e} \BF \in H^{-1/2}_\Psi(\GG_e)$ since $\Lcal_{\Gl, \Gm} \BF=0$ in $D_e$.

Following section~2 in \cite{ACKLM} we define $\Kbb^*: \Hscr^* \to \Hscr^*$ as
\beq\label{Kbbdef}
  \Kbb^* :=
  \begin{bmatrix}
    - \BK_{\GG_i}^* & - \p_{\Gv_i} \BS_{\GG_e} \\
    \p_{\Gv_e} \BS_{\GG_i} & \BK_{\GG_e}^*
  \end{bmatrix}.
\eeq
Then the equation \eqref{raw_matrix_equation} can be rewritten as
\beq
  \left( z_\Gd I + \Kbb^* \right) [{\bf \GF}^\Gd] = \BG, \label{int_equation}
\eeq
where
\beq
  {\bf \GF}^\Gd :=
  \begin{bmatrix}
    \BGvf_i^\Gd \\
    \BGvf_e^\Gd
  \end{bmatrix},
  \quad
  \BG :=
  \begin{bmatrix}
    \p_{\Gv_i} \BF \\
    - \p_{\Gv_e} \BF
  \end{bmatrix}. \nonumber
\eeq
The operator $\Kbb^*$ is the NP operator associated with two circles.

Next, we define $\Sbb$ on $\Hscr^*$ as
\beq\label{Sbbdef}
  \Sbb :=
  \begin{bmatrix}
    \BS_{\GG_i} & \BS_{\GG_e} \\
    \BS_{\GG_i} & \BS_{\GG_e}
  \end{bmatrix}.
\eeq
We emphasize that the $\BS_{\GG_e}$ in the first row is defined to be $\BS_{\GG_e}[\BGvf]|_{\GG_i}$ and $\BS_{\GG_i}$ in the second row is $\BS_{\GG_i}[\BGvf]|_{\GG_e}$. Define, for ${\bf \GF}, {\bf \GY} \in \Hscr^*$,
\beq
  \left( {\bf \GF}, {\bf \GY} \right)_{*} := -\left\la {\bf \GF},  \Sbb [{\bf \GY}] \right\ra, \label{inner_product_Hscr}
\eeq
where $\left\la \cdot, \cdot \right\ra$ denotes the pairing of $L^2(\GG_i)^2 \times L^2(\GG_e)^2$.
Following the same arguments in subsection~3.3 of \cite{ACKLM}, one can show that $( \cdot, \cdot )_{*}$ is actually an inner product on $\Hscr_\GY^*$, and $\Kbb^*$ is a self-adjoint operator with respect to this inner product. Let $\| \cdot \|_*$ be the norm induced by the new inner product. One can show in the same way as in \cite{AJKKY} that this norm is equivalent to the usual $H^{-1/2}$-norm.

Let $\Bu^\Gd$ be the solution to \eqnref{negative_indx_elasticity_problem} given in the form \eqnref{rep_sol}. Using Green's formula, namely,
$$
  \int_{\p \GO} \Bu \cdot \p_\Gv \Bv d\sigma = \int_\GO \Bu \cdot \Lcal_{\Gl, \Gm} \Bv + \int_\GO \Gl \left( \nabla \cdot \Bu \right) \left( \nabla \cdot \Bv \right) + 2 \Gm \hatna \Bu : \hatna \Bv,
$$
one can see that
$$
\int_{D_e \setminus D_i} \Gl | \nabla \cdot (\BS_{\GG_i} [\BGvf_i^\Gd] + \BS_{\GG_e} [\BGvf_e^\Gd] ) |^2 + 2 \Gm | \hatna (\BS_{\GG_i} [\BGvf_i^\Gd] + \BS_{\GG_e} [\BGvf_e^\Gd] ) |^2 = \| \BGF^\Gd \|_*^2.
$$
So, we obtain the following lemma
\begin{lemma}\label{lem:Eudel}
There are constants $C_1$ and $C_2$ such that
\beq
C_1 (\| \BGF^\Gd \|_*^2 -1) \le E(\Bu^\Gd) \le C_2 (\| \BGF^\Gd \|_*^2 +1).
\eeq
\end{lemma}

\section{Asymptotics of NP eigenvalues}\label{sec:NP}

In this section we represent the NP operator as block matrices acting on finite dimensional subspaces, and derive asymptotics of the NP eigenvalues. To do so, we first recall some computations in \cite{AJKKY}.

\subsection{Computations on a disk}

Let $\GG$ be the circle of radius $r_0$. For an integer $m$ let
\beq
\BGvf_m = \begin{bmatrix} \cos m\Go \\
    \sin m\Go \end{bmatrix} \quad\mbox{and}\quad
\widetilde \BGvf_m = \begin{bmatrix} - \sin m\Go \\
    \cos m\Go \end{bmatrix}.
\eeq
The following computations are from \cite{AJKKY}. For $\Bx=(r\cos\Go, r \sin \Go)$,
if $m=1$, then
\beq\label{oneGvf}
- \BS_{\GG} [\BGvf_1](\Bx) =
\begin{cases}
\ds \frac{\Ga_1 - \Ga_2}{2} r \BGvf_1(\Go), & \quad r < r_0,\\
\ds \frac{\Ga_1 - \Ga_2}{2} \frac{r_0^2}{r} \BGvf_1(\Go), & \quad r > r_0.
\end{cases}
\eeq
If $m \ge 2$, then
\beq\label{twoGvf}
- \BS_{\GG} [\BGvf_m](\Bx) =
  \begin{cases}
    \ds \frac{\Ga_1}{2 m} \frac{r^{ m}}{r_0^{m-1}} \BGvf_m(\Go) + \frac{\Ga_2}{2} \left( \frac{r_0^2}{r} - r \right) \frac{r^{m-1}}{r_0^{m-1}} \BGvf_{-m+2}(\Go), & \quad r < r_0, \\
    \ds \frac{\Ga_1}{2 m} \frac{r_0^{m + 1}}{r^{m}} \BGvf_m(\Go), & \quad r > r_0.
  \end{cases}
\eeq
If $m \le -1$, then
\beq\label{negaGvf}
- \BS_{\GG} [\BGvf_m](\Bx) =
  \begin{cases}
    \ds \frac{\Ga_1}{2|m|} \frac{r^{|m|}}{r_0^{|m| - 1}} \BGvf_m(\Go), & \quad r < r_0, \\
    \ds \frac{\Ga_1}{2|m|} \frac{r_0^{|m|+1}}{r^{|m|}} \BGvf_m(\Go) + \frac{\Ga_2}{2} \left( r - \frac{r_0^2}{r} \right) \frac{r_0^{|m|+1}}{r^{|m|+1}} \BGvf_{-m+2}, & \quad r > r_0.
  \end{cases}
\eeq

It is also shown that if $m \ge 2$, then
\beq\label{posiwideGvf}
- \BS_{\GG} [\widetilde\BGvf_m](\Bx) =
  \begin{cases}
  \ds \frac{\Ga_1}{2m} \frac{r^m}{r_0^{m-1}} \widetilde\BGvf_m(\Go) - \frac{\Ga_2}{2} \left( \frac{r_0^2}{r} - r \right) \frac{r^{m-1}}{r_0^{m-1}} \widetilde\BGvf_{-m+2}(\Go), & \quad r < r_0, \\
    \ds \frac{\Ga_1}{2m} \frac{r_0^{m + 1}}{r^{m}}\widetilde\BGvf_m(\Go) , & \quad r > r_0.
  \end{cases}
\eeq
If $m \le -1$, then
\beq\label{negawideGvf}
- \BS_{\GG} [\widetilde\BGvf_m](\Bx) =
  \begin{cases}
      \ds \frac{\Ga_1}{2|m|} \frac{r^{|m|}}{r_0^{|m|-1}} \widetilde\BGvf_m(\Go), & \quad r < r_0, \\
      \ds \frac{\Ga_1}{2|m|} \frac{r_0^{|m|+1}}{r^{|m|}} \widetilde\BGvf_m(\Go) - \frac{\Ga_2}{2} \left( r - \frac{r_0^2}{r} \right) \frac{r_0^{|m|+1}}{r^{|m|+1}} \widetilde\BGvf_{-m+2}(\Go), & \quad r > r_0 .
  \end{cases}
\eeq

As consequences of these computations, eigenvalues of $\BK^*_{\GG}$ are obtained in \cite{AJKKY}:
\beq\label{NPeigen1}
\BK^*_{\GG}[\BGvf_m]= -k_0 \BGvf_m \ (m <0), \quad  \BK^*_{\GG}[\widetilde \BGvf_m]= -k_0 \widetilde\BGvf_m \ (m <0),
\eeq
and
\beq
\BK^*_{\GG}[\BGvf_1]= - \frac{\Gl}{2 (2 \Gm + \Gl)} \BGvf_1, \quad \BK^*_{\GG}[\widetilde \BGvf_m]= k_0 \widetilde\BGvf_m \ (m \ge 2).
\eeq
Even though it is irrelevant to this work, we mention for the completeness's sake that
$$
\BK^*_{\GG}[\widetilde \BGvf_1]= \frac{1}{2} \widetilde\BGvf_1.
$$

\subsection{Block matrix structure of the inner product}

For $n=0, 1,2, \ldots$, we define
\begin{align}
  \BGF_{n, 1} := \begin{bmatrix} \BGvf_{n+1} \\ {\bf 0} \end{bmatrix}, \quad
  \BGF_{n, 2} := \begin{bmatrix} \BGvf_{-n+1} \\ {\bf 0} \end{bmatrix}, \quad
  \BGF_{n, 3} := \begin{bmatrix} {\bf 0} \\ \BGvf_{n+1} \end{bmatrix}, \quad
  \BGF_{n, 4} := \begin{bmatrix} {\bf 0} \\ \BGvf_{-n+1} \end{bmatrix}.
\end{align}
Here ${\bf 0}$ denotes two-dimensional zero vector, and so $\BGF_{n, j}$ are four-dimensional vector-valued functions.
We then define finite dimensional subspaces $\Vscr_n$ of $\Hcal^*$ by
\beq
  \Vscr_n :=
  \mbox{span} \left\{ \BGF_{n, 1}, \BGF_{n, 2}, \BGF_{n, 3}, \BGF_{n, 4} \right\}, \quad n \geq 0.
\eeq
It is worth mentioning that $\Vscr_0$ is of two dimensions and spanned by $\BGF_{0,1}$ and $\BGF_{0,3}$.

Similarly, we define
\begin{align}
  \widetilde \BGF_{n, 1} := \begin{bmatrix} \widetilde \BGvf_{n+1} \\ {\bf 0} \end{bmatrix}, \quad
  \widetilde \BGF_{n, 2} := \begin{bmatrix} \widetilde \BGvf_{-n+1} \\ {\bf 0} \end{bmatrix}, \quad
  \widetilde \BGF_{n, 3} := \begin{bmatrix} {\bf 0} \\ \widetilde \BGvf_{n+1} \end{bmatrix}, \quad
  \widetilde \BGF_{n, 4} := \begin{bmatrix} {\bf 0} \\ \widetilde \BGvf_{-n+1} \end{bmatrix},
\end{align}
and define
\beq
  \widetilde \Vscr_n=
    \mbox{span} \{ \widetilde \BGF_{n, 1}, \widetilde \BGF_{n, 2}, \widetilde \BGF_{n, 3}, \widetilde \BGF_{n, 4} \}, \quad n \geq 0.
\eeq

Using \eqnref{oneGvf}-\eqnref{negawideGvf} in the previous subsection, one can see that
\beq\label{inclusion}
\Sbb(\Vscr_n) \subset \Vscr_n \quad\mbox{and}\quad \Sbb(\widetilde\Vscr_n) \subset \widetilde\Vscr_n \quad\mbox{for each } n.
\eeq
For example, according to the definition \eqnref{Sbbdef}, $\Sbb[\BGF_{n, 1}]$ is given by
$$
\Sbb[\BGF_{n, 1}]=
  \begin{bmatrix}
    \BS_{\GG_i} [\BGvf_{n+1}] |_{\GG_i} \\
    \nm
    \BS_{\GG_i} [\BGvf_{n+1}] |_{\GG_e}
  \end{bmatrix}.
$$
So, using \eqnref{oneGvf}-\eqnref{negaGvf}, one can see that $\Sbb[\BGF_{n, 1}] \in \Vscr_n$.

Moreover, one can derive the block matrix representing $-\Sbb$ on each $\Vscr_n$ or $\widetilde\Vscr_n$. In the following, that $-\Sbb$ is represented by a block matrix $(a_{ij})_{i,j=1}^4$ on $\Vscr_n$ means
$$
-\Sbb[\BGF_{n,i}] = \sum_{j} a_{ij} \BGF_{n,j}, \quad i=1, \ldots, 4.
$$
In fact, they are given as follows. We do not include details of derivations since they are straight-forward, but tedious.
Set $\Gr := r_i / r_e$.
\begin{itemize}
\item[(i)] On $\Vscr_n$, $n \ge 2$,
\beq\label{Sblock}
-\Sbb =
  \begin{bmatrix}
    \frac{\Ga_1}{2 \left( n + 1 \right)} r_i & 0 &  \frac{\Ga_1}{2 \left( n + 1 \right)} r_i \Gr^n & 0 \\
    0 & \frac{\Ga_1}{2 \left( n - 1 \right)} r_i & \frac{\Ga_2}{2}(\frac{r_e^2}{r_i} - r_i) \Gr^n & \frac{\Ga_1}{2 \left( n - 1 \right)} \frac{r_e^2}{r_i} \Gr^n \\
    \frac{\Ga_1}{2 \left( n + 1 \right)} \frac{r_i^2}{r_e} \Gr^n &  \frac{\Ga_2}{2} \left( r_e - \frac{r_i^2}{r_e} \right) \Gr^n &  \frac{\Ga_1}{2 \left( n + 1 \right)} r_e & 0 \\
    0  &  \frac{\Ga_1}{2 \left( n - 1 \right)} r_e \Gr^n  & 0 &  \frac{\Ga_1}{2 \left( n - 1 \right)} r_e
  \end{bmatrix} .
\eeq

\item[(ii)] On $\widetilde{\Vscr}_n$, $n \ge 2$,
\beq\label{wideSblock}
  -\Sbb =
  \begin{bmatrix}
    \frac{\Ga_1}{2 \left( n + 1 \right)} r_i & 0 &  \frac{\Ga_1}{2 \left( n + 1 \right)} r_i \Gr^n & 0 \\
    0 &  \frac{\Ga_1}{2 \left( n - 1 \right)} r_i & -\frac{\Ga_2}{2} \left( \frac{r_e^2}{r_i} - r_i \right) \Gr^n &  \frac{\Ga_1}{2 \left( n - 1 \right)} \frac{r_e^2}{r_i} \Gr^n \\
    \frac{\Ga_1}{2 \left( n + 1 \right)} \frac{r_i^2}{r_e} \Gr^n & -\frac{\Ga_2}{2} \left( r_e - \frac{r_i^2}{r_e} \right) \Gr^n & \frac{\Ga_1}{2 \left( n + 1 \right)} r_e & 0 \\
    0 & \frac{\Ga_1}{2 \left( n - 1 \right)} r_e \Gr^n  & 0 & \frac{\Ga_1}{2 \left( n - 1 \right)} r_e
  \end{bmatrix}.
\eeq
\end{itemize}
We mention that only the cases when $n \ge 2$ are presented above. It is because the cases of large $n$ matter for our analysis.

As a consequence we obtain the following lemma.

\begin{lemma}\label{lem:decomp_Hscr}
\begin{itemize}
\item[(i)] It holds that $\Vscr_m \perp \Vscr_n$ and $\widetilde{\Vscr}_m \perp \widetilde{\Vscr}_n$ if $m\neq n$, and $\Vscr_m \perp \widetilde{\Vscr}_n$ for all $m$ and $n$, with respect to the inner product $( \cdot, \cdot )_{*}$.
\item[(ii)] For $j=1, \ldots, 4$,
\beq\label{phinest1}
( \BGF_{n, j}, \BGF_{n, j} )_{*} \approx \frac{1}{n}, \quad ( \widetilde\BGF_{n, j}, \widetilde\BGF_{n, j} )_{*} \approx \frac{1}{n}  .
\eeq
\item[(iii)] If $i \neq j$, then
\beq\label{phinest2}
      \left| ( \BGF_{n, i}, \BGF_{n, j} )_{*} \right| \lesssim \Gr^n, \quad \left| ( \widetilde\BGF_{n, i}, \widetilde\BGF_{n, j} )_{*} \right| \lesssim \Gr^n .
\eeq
  \end{itemize}
\end{lemma}
In the above and afterwards $( \BGF_{n, j}, \BGF_{n, j} )_{*} \approx \frac{1}{n}$ mean that there are constants $C_1$ and $C_2$ independent of $n$ such that
$$
\frac{C_1}{n} \le \left( \BGF_{n, j}, \BGF_{n, j}\right)_{*} \le \frac{C_2}{n},
$$
and $| ( \BGF_{n, i}, \BGF_{n, j} )_{*} | \lesssim \Gr^n$ means that there is a constant $C$ independent of $n$ such that
$$
\left| ( \BGF_{n, i}, \BGF_{n, j} )_{*} \right| \leq C \Gr^n.
$$

\noindent {\sl Proof of Lemma \ref{lem:decomp_Hscr}}.
The item (i) is an immediate consequence of \eqnref{inclusion}. For (ii) and (iii), we note that $\la \BGF_{n, i}, \BGF_{n, j}\ra_{L^2}=0$ if $i \neq j$, and $\la \BGF_{n, j}, \BGF_{n, j}\ra_{L^2}$ is $2\pi r_i$ if $j=1,2$, and $2\pi r_e$ if $j=3,4$. We then observe that the diagonal elements of the matrices in \eqnref{Sblock} and \eqnref{wideSblock} are of order $n^{-1}$, and off-diagonal elements are of order $\Gr^n$. So, \eqnref{phinest1} and \eqnref{phinest2} follow.
\qed

\subsection{Block matrix structure of the NP operator}

We now represent the NP operator $\Kbb^*$ on $\Vscr_n$ and $\widetilde\Vscr_n$.
According to the definition \eqnref{Kbbdef}, $\Kbb^*[\BGF_{n, 1}]$, for examples, is given by
$$
\Kbb^* [\BGF_{n, 1}] =
  \begin{bmatrix}
    - \BK_{\GG_i}^* [\BGvf_{n+1}] |_{\GG_i} \\
    \nm
    \p_{\Gv_e} \BS_{\GG_i} [\BGvf_{n+1}] |_{\GG_e}
  \end{bmatrix}.
$$
So, using \eqnref{twoGvf} and \eqnref{NPeigen1}, one can compute $\Kbb^* [\BGF_{n, 1}]$. In this way, one can show that
\beq\label{Kinclusion}
\Kbb^*(\Vscr_n) \subset \Vscr_n \quad\mbox{and}\quad \Kbb^*(\widetilde\Vscr_n) \subset \widetilde\Vscr_n \quad\mbox{for each } n.
\eeq
We also obtain the following lemma for the block matrix representation of the NP operator. We again omit the proof since it is straight-forward, but tedious.

\begin{lemma}\label{lem:NPblock}
The NP-operator $\Kbb^*$ admits the following matrix representation:
\begin{itemize}
\item[(i)] On $\Vscr_n$, $n \ge 2$,
\beq\label{J+Mn}
  \Kbb^* = J + M_n,
\eeq
where
\begin{equation*}
  J = k_0
  \begin{bmatrix}
    -1 && 0 && 0 && 0 \\
    \ds 0 && 1 && 0 && 0 \\
    \ds 0 && 0 && 1 && 0 \\
    \ds 0 && 0 && 0 && -1
  \end{bmatrix}
\end{equation*}
and
\begin{equation*}
  M_n = \Gr^n
  \begin{bmatrix}
    0 & 0 & \Gm \Ga_2 & 0 \\
    0 & 0 & - \left( n - 1 \right) \Gm \Ga_2 \left( 1 - \Gr^{-2} \right) & \Gm \Ga_1 \Gr^{-2} \\
    \Gm \Ga_1 \Gr^2 & - \left( n + 1 \right) \Gm \Ga_2 \left( \Gr^2 - 1 \right) & 0 & 0 \\
    0 & \Gm \Ga_2 & 0 & 0
  \end{bmatrix}.
\end{equation*}

\item[(ii)] On $\widetilde{\Vscr_n}$, $n \ge 2$:
\begin{equation*}
  \Kbb^* = J + \widetilde M_n,
\end{equation*}
where
\begin{equation*}
  \widetilde M_n = \Gr^n
  \begin{bmatrix}
    0 & 0 & \Gm \Ga_2 & 0 \\
    0 & 0 & \left( n - 1 \right) \Gm \Ga_2 \left( 1 - \Gr^{- 2} \right) & \Gm \Ga_1 \Gr^{- 2} \\
    \Gm \Ga_1 \Gr^2 & \left( n + 1 \right) \Gm \Ga_2 \left( \Gr^2 - 1 \right) & 0 & 0 \\
    0 & \Gm \Ga_2 & 0 & 0
  \end{bmatrix}.
\end{equation*}
\end{itemize}
\end{lemma}

\subsection{Asymptotics of NP eigenvalues}

It turns out that the exact expression of the eigenvalues of $\Kbb^*$ on $\Vscr_n$, or of the matrix $J+ M_n$ in \eqnref{J+Mn}, is extremely lengthy and complicated. However, for analysis of this paper, it is enough
to have their asymptotic behavior as $n \to \infty$, which we obtain in this subsection using the perturbation theory.

To investigate the asymptotic behavior of the eigenvalues, it is more convenient to express $\Kbb^*= J + M_n$ on $\Vscr_n$ as follows:
\begin{equation*}
  \Kbb^* = P_n + \Gr^n Q,
\end{equation*}
where
\begin{equation*}
  P_n =
  \begin{bmatrix}
    -k_0 & 0 & 0 & 0 \\
    \ds 0 & k_0 & -n \rho^n \mu\alpha_2 (1-\rho^{-2}) & 0 \\
    \ds 0 & -n \rho^n \mu\alpha_2(\rho^2-1) & k_0 & 0 \\
    \ds 0 & 0 & 0 & -k_0
  \end{bmatrix},
\end{equation*}
and
\begin{equation*}
  Q =
  \begin{bmatrix}
    0 & 0 & \Gm \Ga_2 & 0 \\
    0 & 0 & \Gm \Ga_2 \left( 1 - \Gr^{-2} \right) & \Gm \Ga_1 \Gr^{-2} \\
    \Gm \Ga_1 \Gr^2 & - \Gm \Ga_2 \left( \Gr^2 - 1 \right) & 0 & 0 \\
    0 & \Gm \Ga_2 & 0 & 0
  \end{bmatrix}.
\end{equation*}
It is worth emphasizing that $Q$ is independent of $n$.

It is easy to find exact eigenvalues and eigenfunctions of $P_n$. In fact, they are given as follows:
\begin{align}
&\Gl_{n,1}^0 = -k_0, \qquad\qquad\qquad\qquad\qquad \quad \BE_{n,1}^0=[ 1, 0, 0, 0]^T,
\\
&\Gl_{n,2}^0 = +k_0 -n \rho^n \Gm \Ga_2(\Gr-\Gr^{-1}), \qquad \BE_{n,2}^0=[0,1,\Gr,0]^T,  \label{BEn2} \\
&\Gl_{n,3}^0 = +k_0 +n \rho^n \Gm \Ga_2(\Gr-\Gr^{-1}), \qquad \BE_{n,3}^0=[0,-1,\Gr,0]^T,
\\
&\Gl_{n,4}^0 = -k_0, \qquad \qquad\qquad\qquad\quad\quad \quad \BE_{n,4}^0= [ 0, 0, 0, 1]^T.
\end{align}
Here and afterwards, $\BE_{n,j}^0$ written as a vector in $\Rbb^4$ actually represents a four dimensional vector-valued function. For example, $\BE_{n,2}^0$ in \eqnref{BEn2} represents $\BGF_{n,2} + \Gr \BGF_{n,3}$.

Note that, if $n$ is large, the matrix $\Gr^n Q$ becomes a small perturbation matrix. So we can derive asymptotic formula for eigenvalues using  standard arguments of the eigenvalue perturbation theory (see, for example, section XII of \cite{RS}).
We only mention that since eigenvalues $\Gl_2^0$ and $\Gl_3^0$ are simple, we apply non-degenerate perturbation theory. For
the cases of $\Gl_1^0=\Gl_4^0$, we apply the degenerate perturbation theory.

Let $\Gl_{n,j}$ and $\BE_{n,j}$ $(n \ge 1, j=1,2,3,4)$
be the eigenvalues and eigenvectors of $\mathbb{K}^*$, respectively. Then we have
\begin{align}
&\Gl_{n,1} = -k_0 - \rho^{2n} \frac{\mu^2\alpha_1\alpha_2 \rho^4}{k_0(1+\rho^2)}  +O(\rho^{3n}), \qquad \ \BE_{n,1}=\BE_{n,1}^0+O(n\rho^n),
\label{Gln1} \\
&\Gl_{n,2} = +k_0 -n \rho^n \Gm \Ga_2(\Gr-\Gr^{-1})+ O(\rho^n), \qquad \BE_{n,2}=\BE_{n,2}^0+O(\rho^n),
\label{Gln2} \\
&\Gl_{n,3}= +k_0 +n \rho^n \Gm \Ga_2(\Gr-\Gr^{-1})+ O(\rho^n), \qquad \BE_{n,3}=\BE_{n,3}^0+O(\rho^n),
\label{Gln3} \\
&\Gl_{n,4} = -k_0 - \rho^{2n} \frac{\mu^2\alpha_1\alpha_2 \rho^{-2}}{k_0(1+\rho^2)}  +O(\rho^{3n}),\qquad \BE_{n,4}=\BE_{n,4}^0+O(n\rho^n). \label{Gln4}
\end{align}

Similarly, we can also obtain the asymptotic behavior of eigenvalues of $\mathbb{K}^*$ on $\widetilde{\Vscr}_n$ as follows
\begin{align}
&\widetilde{\Gl}_{n,1} = -k_0 - \rho^{2n} \frac{\mu^2\alpha_1\alpha_2 \rho^4}{k_0(1+\rho^2)}  +O(\rho^{3n}),
\label{tildeGln1} \\
&\widetilde\Gl_{n,2}= +k_0 + n \rho^n \Gm \Ga_2(\Gr-\Gr^{-1}) +O(\rho^n),
\\
& \widetilde\Gl_{n,3} = +k_0 - n \rho^n \Gm \Ga_2(\Gr-\Gr^{-1}) + O(\rho^n),
\\
&\widetilde{\Gl}_{n,4} = -k_0 - \rho^{2n} \frac{\mu^2\alpha_1\alpha_2 \rho^{-2}}{k_0(1+\rho^2)}  +O(\rho^{3n}), \label{tildeGln4}
\end{align}
and the corresponding eigenvectors are
\begin{align}
&\widetilde{\BE}_{n,1}=[1, 0, 0, 0]^T+O(n\rho^n), \label{tildeBE1}
\\
&\widetilde{\BE}_{n,2}=[0,1,\rho,0]^T+O(\rho^n),
\\
&\widetilde{\BE}_{n,3}=[0,-1,\rho,0]^T+O(\rho^n),
\\
&\widetilde{\BE}_{n,4}=[0, 0, 0, 1]^T+O(n\rho^n). \label{tildeBE4}
\end{align}
We emphasize that all eigenvalues in \eqnref{Gln1}--\eqnref{tildeGln4} converge to either $k_0$ or $-k_0$.

In view of \eqnref{phinest1}, \eqnref{phinest2}, \eqnref{Gln1}--\eqnref{Gln4}, and \eqnref{tildeBE1}--\eqnref{tildeBE4}, we obtain the following lemma.

\begin{lemma}\label{eigen_inner_estim}
We have
\begin{align}
( \BE_{n, j}, \BE_{n, j} )_{*} & \approx n^{-1}, \quad ( \widetilde\BE_{n, j}, \widetilde\BE_{n, j})_{*} \approx n^{-1} \quad \mbox{for } 1 \leq j \leq 4,  \\
|( \BE_{n, j}, \BGF_{n, j} )_{*}| & \approx  n^{-1}, \quad  |( \widetilde\BE_{n, j}, \widetilde\BGF_{n, j} )_{*}| \approx  n^{-1} \quad \mbox{for } j=1,4, \\
| ( \BE_{n, j}, \BGF_{n, k} )_{*} | & \lesssim \Gr^n, \quad | ( \widetilde\BE_{n, j}, \widetilde\BGF_{n, k} )_{*} | \lesssim \Gr^n \quad \mbox{for }  j=1,4, \, k\neq j, \\
|( \BE_{n, j}, \BGF_{n, k} )_{*}|  & \approx  n^{-1}, \quad |( \widetilde\BE_{n, j}, \widetilde\BGF_{n, k} )_{*}|  \approx  n^{-1} \quad \mbox{for }j=2,3, \,k=2,3, \\
| ( \BE_{n, j}, \BGF_{n, k} )_{*} | & \lesssim \Gr^n, \quad | ( \widetilde\BE_{n, j}, \widetilde\BGF_{n, k} )_{*} | \lesssim \Gr^n \quad \mbox{for }  j=2,3, \,k=1,4.
\end{align}
\end{lemma}

\section{CALR}\label{sec:CALR}

\subsection{Estimates of the solution to the integral equation}

We now look into the integral equation \eqref{int_equation}. We first observe that since $\BG \in \Hscr_\GY^*$, $\BG$ is orthogonal to $\widetilde{\Vscr}_0$ in particular. So $\BG$ is uniquely represented as
$$
  \BG= \sum_{n = 0}^\infty \BG_n + \sum_{n = 1}^\infty \widetilde{\BG}_n, \quad \BG_n \in \Vscr_n, \ \widetilde{\BG}_n \in \widetilde{\Vscr}_n.
$$
It may be helpful to mention that $\BG_1$ has no component of $\BGF_{1,2}$ and $\BGF_{1,4}$ since $\BG$ is perpendicular to constant vectors. The solution $\BGF^\Gd$ to \eqref{int_equation} is given by
\beq
  {\bf \GF}^\delta= \sum_{n = 0}^\infty \BGF_n^\Gd + \sum_{n=1}^\infty \widetilde{\BGF}_n^\Gd, \label{BGFdelta}
\eeq
where $\BGF_n^\Gd \in \Vscr_n$ and $\widetilde{\BGF}_n^\Gd \in \widetilde{\Vscr}_n$ are, respectively, the solutions to
the finite dimensional equations
\beq \label{Int_eqn_mat}
  (z_{\Gd} I + \Kbb^*) \BGF_n^\Gd = \BG_n \ \text{ on } \Vscr_n, \quad (z_{\Gd} I + \Kbb^*) \widetilde{\BGF}_n^\Gd = \widetilde{\BG}_n \ \text{ on }
  \widetilde{\Vscr}_n.
\eeq

Let us consider the first equation in the above.
Since $\BE_{n,j}$, $1 \le j \le 4$, is an orthogonal basis for $\Vscr_n$, we have
$$
\BGF_n^\Gd = \sum_{j=1}^4 \frac{(\BG_n,
\BE_{n,j})_*}{(\BE_{n,j},\BE_{n,j})_* (z_{\Gd} + \Gl_{n,j})}
\BE_{n,j},
$$
and hence
\beq\label{BGFn1}
\| \BGF_n^\Gd \|_*^2 = \sum_{j=1}^4 \frac{|(\BG_n,\BE_{n,j})_*|^2}{\| \BE_{n,j} \|_*^2 |z_{\Gd} + \Gl_{n,j}|^2}.
\eeq

Note that $z_\Gd \to z(c)$ as $\Gd \to 0$ where $z(c)$ is defined in \eqnref{ccond}.
On the other hand, $\Gl_{n,j}$ approaches to either $k_0$ or $-k_0$ as $n \to \infty$. So, if $z(c) \neq \pm k_0$, then $|z_\Gd + \Gl_{n,j}| \ge C$ for some constant $C$ for all sufficiently large $n$ if $\Gd$ is small. So the norm given in \eqnref{BGFn1} does not blow up. So, we assume that $z(c)$ is either $k_0$ or $-k_0$.

Suppose that $z(c)=k_0$. If $j=2,3$, then we infer from \eqnref{Gln2} and \eqnref{Gln3} that $|z_\Gd + \Gl_{n,j}| \ge C$ for some constant $C$ independent of $n$. So, we have
\beq\label{BGFn2}
\sum_{j=2,3} \frac{|(\BG_n,\BE_{n,j})_*|^2}{\| \BE_{n,j} \|_*^2 |z_{\Gd} + \Gl_{n,j}|^2} \lesssim \| \BG_n \|_*^2.
\eeq
If $j=1,4$, then we infer from \eqnref{Gln1} and \eqnref{Gln4} that
\beq\label{BGFn3}
|z_\Gd + \Gl_{n,j}|^2 \approx \Gd^2 + \rho^{4n},
\eeq
and hence
\beq\label{BGFn4}
\sum_{j=1,4} \frac{|(\BG_n,\BE_{n,j})_*|^2}{\| \BE_{n,j} \|_*^2 |z_{\Gd} + \Gl_{n,j}|^2} \approx \sum_{j=1,4} \frac{|(\BG_n,\BE_{n,j})_*|^2}{\| \BE_{n,j} \|_*^2 (\Gd^2 + \rho^{4n})} .
\eeq

Since $\BG_n \in \Vscr_n$ and $\widetilde \BG_n \in \widetilde \Vscr_n$, we have
\beq\label{BGndef}
\BG_n = \sum_{k=1}^4 g_{n,k} \BGF_{n,k}, \quad \widetilde \BG_n = \sum_{k=1}^4 \widetilde g_{n,k} \widetilde \BGF_{n,k},
\eeq
for some constants $g_{n,k}$ and $\widetilde g_{n,k}$. Since $(\BG_{n}, \BE_{n,1})_{*} = \sum_{k=1}^4 g_{n,k}(\BGF_{n,k}, \BE_{n,1})_{*}$, one can see from Lemma \ref{eigen_inner_estim} that
\beq\label{BGFn5}
\frac{|g_{n,1}|^2}{n^2} - \rho^{2n} \sum_{k \neq 1} |g_{n,k}|^2 \lesssim |(\BG_{n}, \BE_{n,1})_{*}|^2 \lesssim \frac{|g_{n,1}|^2}{n^2} + \rho^{2n} \sum_{k \neq 1} |g_{n,k}|^2.
\eeq
Likewise, we have
\beq\label{BGFn6}
\frac{|g_{n,4}|^2}{n^2} - \rho^{2n} \sum_{k \neq 4} |g_{n,k}|^2 \lesssim |(\BG_{n}, \BE_{n,4})_{*}|^2 \lesssim \frac{|g_{n,4}|^2}{n^2} + \rho^{2n} \sum_{k \neq 4} |g_{n,k}|^2.
\eeq
It then follows from Lemma \ref{eigen_inner_estim} and \eqnref{BGFn4}--\eqnref{BGFn6} that
\begin{align}
&\frac{n^{-1} (|g_{n,1}|^2+ |g_{n,4}|^2) - n\rho^{2n} (|g_{n,2}|^2+ |g_{n,3}|^2) }{\Gd^2 + \rho^{4n}} \lesssim \nonumber \\
& \sum_{j=1,4} \frac{|(\BG_n,\BE_{n,j})_*|^2}{\| \BE_{n,j} \|_*^2 |z_{\Gd} + \Gl_{n,j}|^2} \lesssim \frac{n^{-1} (|g_{n,1}|^2+ |g_{n,4}|^2) + n\rho^{2n} (|g_{n,2}|^2+ |g_{n,3}|^2) }{\Gd^2 + \rho^{4n}}. \label{BGFn7}
\end{align}
So, we obtain from \eqnref{BGFn1}, \eqnref{BGFn2}, and \eqnref{BGFn7} that
\begin{align}
&\frac{n^{-1} (|g_{n,1}|^2+ |g_{n,4}|^2) - n\rho^{2n} (|g_{n,2}|^2+ |g_{n,3}|^2) }{\Gd^2 + \rho^{4n}} - \| \BG_n \|_*^2 \lesssim \nonumber \\
& \| \BGF_n^\Gd \|_*^2 \lesssim
\frac{n^{-1} (|g_{n,1}|^2+ |g_{n,4}|^2) + n\rho^{2n} (|g_{n,2}|^2+ |g_{n,3}|^2) }{\Gd^2 + \rho^{4n}} + \| \BG_n \|_*^2 . \label{BGFn8}
\end{align}

One can estimate $\| \widetilde\BGF_n^\Gd \|_*^2$ in a similar way to obtain
\begin{align}
&\frac{n^{-1} (|\widetilde g_{n,1}|^2+ |\widetilde g_{n,4}|^2) - n\rho^{2n} (|\widetilde g_{n,2}|^2+ |\widetilde g_{n,3}|^2) }{\Gd^2 + \rho^{4n}} - \| \widetilde \BG_n \|_*^2 \lesssim \nonumber \\
& \| \widetilde \BGF_n^\Gd \|_*^2 \lesssim
\frac{n^{-1} (|\widetilde g_{n,1}|^2+ |\widetilde g_{n,4}|^2) + n\rho^{2n} (|\widetilde g_{n,2}|^2+ |\widetilde g_{n,3}|^2) }{\Gd^2 + \rho^{4n}} + \| \widetilde \BG_n \|_*^2 . \label{BGFn9}
\end{align}

Let, for ease of notation,
\beq\label{Indef}
I_n:= |g_{n,1}|^2+ |g_{n,4}|^2 +|\widetilde g_{n,1}|^2+ |\widetilde g_{n,4}|^2
\eeq
and
\beq\label{IIndef}
II_n:= |g_{n,2}|^2+ |g_{n,3}|^2 +|\widetilde g_{n,2}|^2+ |\widetilde g_{n,3}|^2.
\eeq
Since
$$
\sum_{n=0}^\infty \| \BG_n \|_*^2 + \sum_{n=1}^\infty \| \widetilde \BG_n \|_*^2 = \| \BG \|_*^2,
$$
we infer by summing \eqnref{BGFn8} and \eqnref{BGFn9} over $n$ that there are constants $C_1$ and $C_2$ independent of $\Gd$ such that
\beq\label{BGFn10}
C_1 \left( \sum_{n=1}^\infty \frac{n^{-1} I_n - n\rho^{2n} II_n}{\Gd^2 + \rho^{4n}} - 1 \right) \le \| \BGF^\Gd \|_*^2 \le C_2 \left( \sum_{n=1}^\infty \frac{n^{-1} I_n + n\rho^{2n} II_n}{\Gd^2 + \rho^{4n}} + 1 \right).
\eeq

If $z(c)=-k_0$, then we obtain in a similar way that
\beq\label{BGFn11}
C_1 \left( \sum_{n=1}^\infty \frac{n^{-1} II_n - n\rho^{2n} I_n}{\Gd^2 + n^2\rho^{2n}} - 1 \right) \le \| \BGF^\Gd \|_*^2 \le C_2 \left( \sum_{n=1}^\infty \frac{n^{-1} II_n + n\rho^{2n} I_n}{\Gd^2 + n^2\rho^{2n}} + 1 \right).
\eeq
It is worth emphasizing that the quantity $n(\Gd^2 + n^2\rho^{2n})$ in the denominator is different from that in \eqnref{BGFn10}. This discrepancy is caused by the different asymptotic behaviors of eigenvalues near $k_0$ and $-k_0$ as shown in \eqnref{Gln1}--\eqnref{tildeGln4}.

In summary, we obtain the following proposition from Lemma \ref{lem:Eudel}, \eqnref{BGFn10} and \eqnref{BGFn11}.

\begin{prop}\label{prop:Eudel}
There are constants $C_1$ and $C_2$ independent of $\Gd$ such that
\begin{itemize}
\item[(i)] if $z(c)=k_0$, then
\beq\label{Eudelone}
C_1 \left( \sum_{n=1}^\infty \frac{n^{-1} I_n - n\rho^{2n} II_n}{\Gd^2 + \rho^{4n}} - 1 \right) \le E(\Bu^\Gd) \le C_2 \left( \sum_{n=1}^\infty \frac{n^{-1} I_n + n\rho^{2n} II_n}{\Gd^2 + \rho^{4n}} + 1 \right),
\eeq
\item[(ii)] if $z(c)=-k_0$, then
\beq\label{Eudeltwo}
C_1 \left( \sum_{n=1}^\infty \frac{n^{-1} II_n - n\rho^{2n} I_n}{\Gd^2 + n^2\rho^{2n}} - 1 \right) \le E(\Bu^\Gd) \le C_2 \left( \sum_{n=1}^\infty \frac{n^{-1} II_n + n\rho^{2n} I_n}{\Gd^2 + n^2\rho^{2n}} + 1 \right).
\eeq
\end{itemize}
\end{prop}

\subsection{Resonance by dipole sources}\label{sec:dipole}

In this section we assume that the source function $\Bf$ in \eqnref{negative_indx_elasticity_problem} is given by the dipole-type function, namely,
\beq\label{dipole}
  \Bf (\Bx) = \Bb^T \nabla_{\Bx} \left(
    \begin{bmatrix}
      \Gd_{\Bz}(\Bx) & 0 \\
      0 & \Gd_{\Bz}(\Bx)
    \end{bmatrix} \Ba \right),
\eeq
where $\Bz = (z_1, z_2)^T$ is the location of the dipole outside $D_e$, and $\Ba = (a_1, a_2)^T$, $\Bb = (b_1, b_2)^T$ are constant vectors. With this source function we obtain the following theorem using Proposition \ref{prop:Eudel}. In what follows, the notation $A \sim B$ ($A$ and $B$ are two quantities depending on $\Gd$) indicates that there are constants $C_1$ and $C_2$ independent of $\Gd$ such that $C_1 \le A/B \le C_2$.

\begin{theorem} \label{thm:res+}
Let $\Bf$ be given by \eqref{dipole}, and let $\rho_\Bz = r_e/|\Bz|$. It holds that
\begin{itemize}
\item[(i)] if $z(c) = k_0$, then
\beq\label{itemone}
E(\Bu^\Gd) \sim \left| \ln \Gd \right|^3 \Gd^{\frac{\ln \rho_\Bz}{\ln \rho} - 2} \quad \text{as } \Gd \to 0,
\eeq

\item[(ii)] if $z(c) = -k_0$, then
\beq\label{itemtwo}
E(\Bu^\Gd) \sim \left| \ln \Gd \right|^{-1-\frac{2\ln\rho_\Bz}{\ln \rho}} \Gd^{\frac{2\ln \rho_\Bz}{\ln \rho} - 2} \quad \text{as } \Gd \to 0.
\eeq
\end{itemize}
\end{theorem}

As an immediate consequence we obtain the following corollary.

\begin{cor}\label{cor:res+}
Suppose that $\Bf$ is given by \eqref{dipole}.
\begin{itemize}
\item[(i)] If $z(c) = k_0$, let $r_*:=r_e^2 / r_i$. As $\Gd \to 0$,
\beq
\Gd E(\Bu^\Gd) \to
\begin{cases}
\infty \quad &\mbox{if } r_e < |\Bz| \le r_* , \\
0 \quad &\mbox{if } |\Bz| > r_*.
\end{cases}
\eeq
\item[(ii)] If $z(c) = -k_0$, let $r_{**}:=\sqrt{r_e^3 / r_i}$. As $\Gd \to 0$,
\beq
\Gd E(\Bu^\Gd) \to
\begin{cases}
\infty \quad &\mbox{if } r_e < |\Bz| < r_{**} , \\
0 \quad &\mbox{if } |\Bz| \ge r_{**}.
\end{cases}
\eeq
\end{itemize}
\end{cor}

\noindent{\sl Proof of Theorem \ref{thm:res+}}.
Since $\Bf$ is assumed to be of the form \eqref{dipole}, $\BF$ defined by \eqnref{Fdef} is given by
\beq
  \BF(\Bx) = \Bb^T \nabla_{\Bx} \left( \BGG(\Bx - \Bz) \Ba \right) = - \Bb^T \nabla_\Bz \left( \BGG(\Bx - \Bz) \Ba \right).
\eeq
Let $g_{m, j}$ and $\widetilde g_{m, j}$, $j = 1,2,3,4$, be the coefficients of $\BG= (\p_{\Gv_i} \BF , - \p_{\Gv_e} \BF)^T$ as
defined by \eqnref{BGndef}. It turns out that $g_{n, j}$ and $\widetilde g_{n, j}$ satisfy
\beq\label{gnone}
  g_{n, 1}, \widetilde g_{n, 1} \approx \left( \frac{r_i}{|\Bz|} \right)^n, \quad
  g_{n, 2}, \widetilde g_{n, 2} \approx n^2 \left( \frac{r_i}{|\Bz|} \right)^n,
\eeq
and
\beq\label{gntwo}
  g_{n, 3}, \widetilde g_{n, 3} \approx \left( \frac{r_e}{|\Bz|} \right)^n, \quad
  g_{n, 4}, \widetilde g_{n, 4} \approx n^2 \left( \frac{r_e}{|\Bz|} \right)^n.
\eeq
We include proofs of these estimates in Appendix \ref{proofs}.

Suppose that $z(c)=k_0$. Then we infer from \eqnref{Indef}, \eqnref{IIndef}, \eqnref{gnone} and \eqnref{gntwo} that
\beq
I_n \approx n^4 \Gr_\Bz^{2n} \quad\mbox{and}\quad II_n \approx \Gr_\Bz^{2n}.
\eeq
It then follows from \eqnref{Eudelone} that
\begin{equation*}
E(\Bu^\Gd) \sim \sum_{n = 1}^\infty \frac{n^3 \Gr_\Bz^{2 n}}{\Gd^2 + \Gr^{4 n}}.
\end{equation*}

For small $\Gd$, let $N$ be a positive integer such that
\beq\label{Ndef}
  \Gr^{2 N} \leq \Gd \leq \Gr^{2 \left( N - 1 \right)}.
\eeq
Note that $N \sim |\ln \Gd|$.
Then we have
\begin{equation*}
E(\Bu^\Gd) \sim \sum_{n < N} + \sum_{n \geq N} \sim \sum_{n < N} \frac{n^3 \Gr_\Bz^{2 n}}{\Gr^{4 n}} + \sum_{n > N} \frac{1}{\Gd^2} {n^3 \Gr_\Bz^{2 n}} =: S_1^\Gd + S_2^\Gd.
\end{equation*}
Since
$$
S_1^\Gd \sim \int_1^N y^3 e^{2 y \ln (\Gr_\Bz / \Gr^2) } \, dy ,
$$
we have
\begin{equation*}
S_1^\Gd \sim N^3 e^{2 N \ln (\Gr_\Bz / \Gr^2)} \sim \left| \ln \Gd \right|^3 \Gd^{\frac{ \ln \Gr_\Bz}{\ln \Gr}-2}.
\end{equation*}
Likewise, since
\begin{equation*}
S_2^\Gd \sim \Gd^{- 2} \int_N^\infty y^3 e^{2 y \ln \Gr_\Bz} \, dy ,
\end{equation*}
we have
\begin{equation*}
S_2^\Gd \sim \Gd^{-2} N^3 e^{2N \ln \Gr_\Bz}\sim \left| \ln \Gd \right|^3 \Gd^{\frac{ \ln \Gr_\Bz}{\ln \Gr}-2} .
\end{equation*}
This proves \eqnref{itemtwo}

If $z(c)=-k_0$, we have
\beq
E(\Bu^\Gd) \sim \sum_{n = 1}^\infty \frac{\Gr_\Bz^{2 n}}{n(\Gd^2 + n^2 \rho^{2 n})}.
\eeq
In this case we choose $N$, instead of \eqnref{Ndef}, so that 
$$
N\rho^N \leq \delta \leq (N-1)\rho^{N-1}.
$$
Then we have $\delta/N\sim e^{N\ln\rho}$ and $N\sim |\ln\delta|$. Then, in a way similar to above, one can show \eqnref{itemtwo}.
\qed

\subsection{Cloaking due to ALR}

Recall that
\begin{equation*}
  \Bu^\Gd = \BF + \BS_{\GG_i} [\BGvf_i^\Gd] + \BS_{\GG_e} [\BGvf_e^\Gd].
\end{equation*}
In this section, we show that
\beq\label{V_bdd}
  \big| \Bu^\Gd(\Bx) \big| \leq C, \quad \left\vert \Bx \right\vert = r > r_0,
\eeq
for some $r_0 > r_e$.
Note that the above boundedness means that CALR happens, since we have already proved the energy $E^\Gd$ blows up when a point source located inside the crical radius.

For simplicity, we only consider $\BS_{\GG_e} [\BGvf_e^\Gd]$.
The potential $\BS_{\GG_i} [\BGvf_i^\Gd]$ can be estimated in a similar way.
Let us write the solutions of the equations \eqref{Int_eqn_mat} using $\BGF_{n, j}$ and $\widetilde \BGF_{n, j}$.
\beq
\BGF_n^\Gd = \sum_{j=0}^4 \frac{g_{n,j}}{z_\Gd + \Gl_{n,j}} \BGF_{n,j}, \quad \widetilde\BGF_n^\Gd = \sum_{j=1}^4 \frac{\widetilde g_{n,j}}{z_\Gd + \widetilde \Gl_{n,j}} \widetilde \BGF_{n,j}, \nonumber
\eeq
where $g_{n,j}$ and $\widetilde g_{n,j}$ are given by \eqref{BGndef}.
Then, from \eqref{BGFdelta}, we have
\begin{align*}
  \BGvf_e^\Gd & = \sum_{n=0}^\infty \left\{ \frac{g_{n,3} \BGvf_{n+1}}{z_\Gd + \Gl_{n,3}} + \frac{\widetilde g_{n,3} \widetilde \BGvf_{n+1}}{z_\Gd + \widetilde \Gl_{n,3}} \right\} + \sum_{n=1}^\infty \left\{ \frac{g_{n,4} \BGvf_{-n+1}}{z_\Gd + \Gl_{n,4}} + \frac{\widetilde g_{n,4} \widetilde \BGvf_{-n+1}}{z_\Gd + \widetilde \Gl_{n,4}} \right\} \\
  & =: \BGvf_{e,3}^\Gd + \widetilde \BGvf_{e,3}^\Gd + \BGvf_{e,4}^\Gd + \widetilde \BGvf_{e,4}^\Gd.
\end{align*}
Hence, we have
\beq \label{BSBGvf_e}
\BS_{\GG_e}[\BGvf_e^\Gd] = \BS_{\GG_e}[\BGvf_{e,3}^\Gd] + \BS_{\GG_e}[\widetilde \BGvf_{e,3}^\Gd] + \BS_{\GG_e}[\BGvf_{e,4}^\Gd] + \BS_{\GG_e}[\widetilde \BGvf_{e,4}^\Gd].
\eeq

Suppose $z(c) = k_0$.
Let us estimate $\BS_{\GG_e}[\BGvf_{e,4}^\Gd]$ in \eqref{BSBGvf_e}.
The other terms can be estimated in the same manner.
From \eqref{negaGvf}, we have
\beq\label{1000}
|\BS_{\GG_e} [\BGvf_{e, 4}^\Gd] (\Bx)| \lesssim \sum_{n=1}^\infty \frac{|g_{n,4}|}{\Gd + \Gr^{2n}} \frac{r_e^n}{r^{n-1}} \lesssim \sum_{n=1}^\infty n^2 \left(\frac{r_e^2}{r_i}\right)^{2n} \left(\frac{1}{r |\Bz|}\right)^{n-1}.
\eeq
Therefore, for $|\Bx| = r > r_e^3 / r_i^2$, we have $|\BS_{\GG_e}[\BGvf_{e,4}^\Gd](\Bx)| < C$ for some $C > 0$.

Next, suppose $z(c) = - k_0$.
For the same reason, we only estimate $\BS_{\GG_e}[\BGvf_{e,3}^\Gd]$.
From \eqref{twoGvf}, we have
\beq\label{1001}
|\BS_{\GG_e} [\BGvf_{e, 3}^\Gd] (\Bx)| \lesssim \sum_{n=1}^\infty \frac{|g_{n,3}|}{\Gd + n \Gr^n} \frac{r_e^n}{r^{n-1}} \lesssim \sum_{n=1}^\infty \frac{1}{n} \left(\frac{r_e^3}{r_i}\right)^n \left(\frac{1}{r |\Bz|}\right)^{n-1}.
\eeq
Therefore, for $|\Bx| = r > r_e^2 / r_i$, we have $|\BS_{\GG_e}[\BGvf_{e,3}^\Gd](\Bx)| < C$ for some $C > 0$.

Since $\BF(\Bx) \to {\bf 0}$ as $|\Bx| \to \infty$, \eqnref{1000} and \eqnref{1001} together with Theorem~\ref{thm:res+} yield the following theorem.
\begin{theorem}
  Suppose that $\Bf$ be given by \eqref{dipole}.
  \begin{itemize}
  \item[(i)] If $z(c) = k_0$, then CALR occurs with the critical radius $r_* = r_e^2 / r_i$; more precisely, $E^\Gd=\Gd E(\Bu^\Gd) \to \infty$ as $\delta \to 0$ if $r_e < |\Bz| \le r_*$, and there is a constant $C$ such that
    \begin{equation*}
      \big\vert \Bu^\Gd(\Bx) \big\vert < C \quad \text{for } \left\vert \Bx \right\vert = r > r_e^3 / r_i^2 \quad \text{as } \Gd \to 0.
    \end{equation*}
  \item[(ii)] If $z(c) = - k_0$, then CALR occurs with the critical radius $r_{**} = \sqrt{r_e^3 / r_i}$; more precisely, $E^\Gd \to \infty$ as $\Gd \to 0$ if $r_e < |\Bz| < r_{**}$, and there is a constant $C$ such that
    \begin{equation*}
      \big\vert \Bu^\Gd(\Bx) \big\vert < C \quad \text{for } \left\vert \Bx \right\vert = r > r_e^2 / r_i \quad \text{as } \Gd \to 0.
    \end{equation*}
  \end{itemize}
\end{theorem}

\appendix

\section{Proofs of \eqnref{gnone} and \eqnref{gntwo}} \label{proofs}

Here we use the complex representation of the displacement vectors. So we identify $\Ba, \Bb, \Bz, \Bx$ with complex numbers $a=a_1+ia_2$, $b=b_1+ib_2$, $z=z_1+iz_2$ and $x=x_1+ix_2$, respectively, where $\Ba, \Bb, \Bz$ are vectors appearing in the definition of
the dipole source given in \eqnref{dipole}. Let
$$
\Gk = \frac{\Gl + 3 \Gm}{\Gl + \Gm}.
$$
We have
\beq \label{Gamma_a_formula}
    2 \Gm \left[ \left(\BGG(\Bx - \Bz) \mathbf{a}\right)_1 + i \left(\BGG(\Bx - \Bz) \mathbf{a}\right)_2 \right]= \Gk \Gf(x) - x \ol{\Gf'(x)}-\ol{\Gy(x)},
\eeq
where
  \begin{equation*}
    \Gf(x) = \frac{a}{2 \pi \left( \Gk + 1 \right)} \ln(x - z),
  \end{equation*}
  and
  \begin{equation*}
    \Gy(x) = - \frac{\Gk \ol{a}}{2 \pi \left( \Gk + 1 \right)} \ln(x - z) - \frac{a}{2 \pi \left( \Gk + 1 \right)} \left( \frac{\ol{z}}{x - z} - \frac{\ol{a}}{a} \right).
  \end{equation*}
We emphasize that since $z$ lies outside $\ol{D_e}$, $\ln(x - z)$ is well-defined as a function of $x \in \ol{D_e}$.

By a straightforward computation using the complex potential representation of $\BGG(\Bx - \Bz) \Ba$, we can obtain
\begin{align}
  & \left[\left(\p_\Gv \BGG(\Bx - \Bz) \Ba\right)_1 + i \left(\p_\Gv \BGG(\Bx - \Bz) \Ba\right)_2\right] \big|_{\GG_e} \\
  & = \sum_{m = 0}^\infty \Re\left\{ G_m(z, \ol{z}) \right\} e^{i \left( m + 1 \right) \Go} + \Im\left\{ G_m(z, \ol{z}) \right\} i e^{i \left( m + 1 \right) \Go} \nonumber \\
  & \qquad + \sum_{m = 2}^\infty \Re\left\{ H_m(z, \ol{z}) \right\} e^{- i \left( m - 1 \right) \Go} + \Im\left\{ H_m(z, \ol{z}) \right\} i e^{- i \left( m - 1 \right) \Go},
\end{align}
where
\begin{equation*}
  G_m(z, \ol{z}) = - \frac{1}{2 \pi (\Gk + 1)}\left( \frac{a}{z} + \frac{\ol a}{\ol z}\Gd_{m 0} \right) \frac{r_e^m}{z^m},
\end{equation*}
and
\begin{equation*}
  H_m(z, \ol{z}) = \frac{1}{2 \pi (\Gk + 1)}\left[(m-1) \frac{\ol a}{\ol z} - \frac{\ol z}{r_e^2} \left(\Gk a + \frac{\ol a z}{\ol z} (m-1)\right) \right] \frac{r_e^m}{\ol{z}^m}.
\end{equation*}
Since
\begin{equation*}
  \p_\Gv \BF \big|_{\GG_e} = - \Bb^T \nabla_{\Bz} \left( \p_\Gv \BGG(\Bx - \Bz) \Ba \big|_{\GG_e} \right),
\end{equation*}
then we have
\begin{equation*}
  g_{m, 3} = - \Bb^T \nabla_{\Bz} \Re\left\{ G_m(z, \ol{z}) \right\}, \qquad \widetilde g_{m, 3} = - \Bb^T \nabla_{\Bz} \Im\left\{ G_m(z, \ol{z}) \right\},
\end{equation*}
and
\begin{equation*}
  g_{m, 4} = - \Bb^T \nabla_{\Bz} \Re\left\{ H_m(z, \ol{z}) \right\}, \qquad \widetilde g_{m, 4} = - \Bb^T \nabla_{\Bz} \Im\left\{ H_m(z, \ol{z}) \right\}.
\end{equation*}
Estimates in \eqnref{gntwo} follows from these.

Estimates in \eqnref{gnone} can be proved similarly.

\end{document}